\documentclass[envcountsame,envcountresetchap]{svmult}

\usepackage{amsmath,amsfonts,amssymb,amscd,verbatim,delarray,verbatim}
\newcommand{\C}{{\mathbb C}}
\newcommand{\F}{{\mathbb F}}

\newcommand{\Q}{{\mathbb Q}}
\newcommand{\Z}{{\mathbb Z}}

\DeclareMathOperator{\M}{M} 
\DeclareMathOperator{\Out}{Out}
\DeclareMathOperator{\Aut}{Aut} 
 
\DeclareMathOperator{\Inndiag}{Inndiag}
\DeclareMathOperator{\Outdiag}{Outdiag}
\DeclareMathOperator{\Res}{Res}
\DeclareMathOperator{\order}{order}

\smartqed

\input xy
\xyoption{all}

\spnewtheorem{thm}{Theorem}[section]{\bf}{\it}
\spnewtheorem{conj}[thm]{Conjecture}{\bf}{\it}
\spnewtheorem{prop}[thm]{Proposition}{\bf}{\it}
\spnewtheorem{defn}[thm]{Definition}{\bf}{\it}
\spnewtheorem{lemm}[thm]{Lemma}{\bf}{\it}
\spnewtheorem{corol}[thm]{Corollary}{\bf}{\it}
\spnewtheorem{assertionp}[thm]{Assertion}{\bf}{\it}

\spnewtheorem{conclusionsp}[thm]{Conclusion}{\bf}{\rm}
\spnewtheorem{assumptionp}[thm]{Assumption}{\bf}{\rm}
\spnewtheorem{exs}[thm]{Examples}{\bf}{\rm}
\spnewtheorem{exam}[thm]{Example}{\bf}{\rm}
\spnewtheorem{ques}[thm]{Question}{\bf}{\rm}
\spnewtheorem{rem}[thm]{Remark}{\bf}{\rm}
\spnewtheorem{rmks}[thm]{Remarks}{\bf}{\rm}
\spnewtheorem{nota}[thm]{Notations}{\bf}{\rm}
\spnewtheorem{disc}[thm]{Discussion}{\bf}{\rm}
\spnewtheorem{chal}[thm]{Challenge}{\bf}{\rm}

\def\F{{\mathbb{F}}}

\def\Q{{\mathbb{Q}}}

\def\Gm{{{\mathbb G}_{\textrm{m}}}}

\def\cO{{\mathcal O}}

\def\Lt{{\widetilde L}}
\def\Gt{{\widetilde G}}

    \def\Hom{\textrm{Hom}}
           \def\int{\textrm{int}}

\def\Br{\textrm{Br}}

\begin{document}
\setcounter{page}{1}

\title*{The Bogomolov multiplier of finite simple groups}
\titlerunning{The Bogomolov multiplier of finite simple groups}
\author{Boris Kunyavski\u\i}
\authorrunning{B.  Kunyavski\u\i}
\institute{Department of
Mathematics\\
 Bar-Ilan University\\
 52900 Ramat Gan, ISRAEL\\
\texttt{kunyav@macs.biu.ac.il}
}

\maketitle


\begin{abstract}
{}
The subgroup of the Schur multiplier of a finite group $G$
consisting of all cohomology classes whose restriction to any
abelian subgroup of $G$ is zero is called the Bogomolov multiplier
of $G$. We prove that if $G$ is quasisimple or almost simple, its
Bogomolov multiplier is trivial except for the case of certain
covers of $PSL(3,4)$.
\end{abstract}

\section*{Introduction} \label{sec:intro}

A common method for proving that a given algebraic variety $X$ over
a field $k$ is not rational is as follows. We consider some easily
computable object (usually of algebraic nature), which can be
defined functorially on a sufficiently large class of algebraic
varieties and is known to be preserved under birational
transformations ({\it birational invariant}). We calculate its value
for $X$ (or for some $Y$ birationally equivalent to $X$). If this
value is not trivial, i.e., does not coincide with the value of this
birational invariant on the affine or projective space, $X$ is not
rational.

The Brauer group $\Br (X)=H_{\textrm{\'et}}^2(X,\Gm )$, whose
birational invariance in the class of smooth projective varieties
has been established by Grothendieck, turned out to be a very
convenient tool (the Artin--Mumford counter-example to L\"uroth's
problem, based on using this invariant, confirms its power).
Moreover, even if $X$ is not projective, this invariant can be
useful: embed $X$ as an open subset into a smooth projective variety
$Y$ (if the ground field is of characteristic zero, this is always
possible by Hironaka) and compute $\Br (Y)$. If the latter group is
not zero, $Y$ cannot be birational to $\mathbb P^n$, and thus $X$ is
not rational. Note that $\Br (Y)$ depends only on $X$ (and not on
the choice of a smooth projective model $Y$); it is called the {\it
unramified Brauer group of $X$} and denoted by $\Br_{\textrm{nr}}
(X)$. (The reader interested in historical perspective and geometric
context, including more general invariants arising from
higher-dimensional cohomology, is referred to \cite{Sh}, \cite{CTS},
\cite[6.6, 6.7]{GS}, \cite{Bo07}.)

In concrete cases, it may be difficult to construct $Y$ explicitly,
and thus it is desirable to express $\Br_{\textrm{nr}} (X)$ in {\it
intrinsic} terms, i.e., get a formula not depending on $Y$. This
approach was realized by Bogomolov \cite{Bo87} in the case of the
quotient variety $X=V/G$ where $V$ stands for a faithful linear
representation of a linear algebraic group $G$ over $\C$. It turns
out that in this case $\Br_{\textrm{nr}} (X)$ depends solely on $G$
(but not on $V$).

In the present paper we focus on the case where $G$ is a {\it
finite} group. The birational invariant $\Br_{\textrm{nr}} (V/G)$
has been used by Saltman to give a negative answer to Noether's
problem \cite{Sa}. In \cite{Bo87} Bogomolov established an explicit
formula for $\Br_{\textrm{nr}} (V/G)$ in terms of $G$: this group is
isomorphic to $B_0(G)$, the subgroup of the Schur multiplier $\M
(G):=H^2(G, \Q/\Z)$ consisting of all cohomology classes whose
restriction to any abelian subgroup of $G$ is zero. We call $B_0(G)$
the {\it Bogomolov multiplier} of $G$.

In \cite{Sa}, \cite{Bo87} one can find examples of groups $G$ with
nonzero $B_0(G)$ (they are all $p$-groups of small nilpotency
class). In contrast, in \cite{Bo92} Bogomolov conjectured that
$B_0(G)=0$ when $G$ is a finite {\it simple} group. In \cite{BMP} it
was proved that $B_0(G)=0$ when $G$ is of Lie type $A_n$. In the
present paper we prove Bogomolov's conjecture in full generality.

\section{Results} \label{results}

We maintain the notation of the introduction and assume throughout
the paper that $G$ is a {\it finite} group.

We say that $G$ is {\it quasisimple} if $G$ is perfect and its
quotient by the centre $L=G/Z$ is a nonabelian simple group. We say
that $G$ is {\it almost simple} if for some nonabelian simple group
$L$ we have $L\subseteq G\subseteq \Aut L$. Our first observation is

\begin{thm} \label{quasi}
If $G$ is a finite quasisimple group other than a $4$- or $12$-cover
of $PSL(3,4)$, then $B_0(G)=0$.
\end{thm}

\begin{corol} \label{cor:Bogom}
If $G$ is a finite simple group, then $B_0(G)=0$.
\end{corol}

This corollary proves Bogomolov's conjecture.

From Corollary \ref{cor:Bogom} we deduce the following

\begin{thm} \label{almost}
If $G$ is a finite almost simple group, then $B_0(G)=0$.
\end{thm}

\begin{rem}
Following \cite[Ch.~2, \S\S~6,7]{GL}, we call quasisimple groups $Q$
(as in Theorem \ref{quasi}) and almost simple groups $A$ (as in
Theorem \ref{almost}), as well as the extensions of $A$ by $Q$, {\it
decorations} of finite simple groups. It is most likely that one can
complete the picture given in the above theorems, allowing both
perfect central extensions and outer automorphisms, by deducing from
Theorems \ref{quasi} and \ref{almost} that $B_0(G)=0$ for all {\it
nearly simple} groups $G$ (see the definition in Section
\ref{sec:fin} below) excluding the groups related to the above
listed exceptional cases. In particular, this statement holds true
for all finite ``reductive'' groups such as the general linear group
$GL(n,q)$, the general unitary group $GU(n,q)$, and the like.
\end{rem}

Our notation is standard and mostly follows \cite{GLS}. Throughout
below ``simple group'' means ``finite nonabelian simple group''. Our
proofs heavily rely on the classification of such groups.

\section{Preliminaries} \label{sec:pre}

In order to make the exposition as self-contained as possible,
in this section we collect the group-theoretic information
needed in the proofs. All groups are assumed finite (although
some of the notions discussed below can be defined for infinite
groups as well).

\subsection{Schur multiplier}  \label{sec:Schur}

The material below (and much more details) can be found in \cite{Ka}.

The group $\M (G):=H^2(G,\Q/\Z)$, where $G$ acts on $\Q/\Z$ trivially,
is called the Schur multiplier of $G$.
It can be identified with the kernel of some central extension
$$
1\to \M (G) \to \widetilde G\to G \to 1.
$$
The covering group $\widetilde G$ is defined uniquely up to
isomorphism provided $G$ is perfect (i.e., coincides with its derived
subgroup $[G,G]$).

We will need to compute $\M (G)$ in the case where $G$ is a
semidirect product of a normal subgroup $N$  and a subgroup $H$. If
$A$ is an abelian group on which $G$ acts trivially, the restriction
map $\Res_H\colon H^2(G,A)\to H^2(H,A)$ gives rise to a split exact
sequence \cite[Prop.~1.6.1]{Ka}
$$
1\to K \to H^2(G,A)\to H^2(H,A)\to 1.
$$
The kernel $K$ can be computed from the exact sequence
\cite[Th.~1.6.5(ii)]{Ka}
$$
1\to H^1(H,\Hom (N,A))\to K \stackrel{\Res_N}{\to} H^2(N,A)^H\to H^2(H,\Hom (N,A)).
$$
If $N$ is perfect and $A=\Q /\Z$, we have $\Hom (N,A)=1$ and thus
\cite[Lemma 16.3.3]{Ka}
\begin{equation}
\M(G)\cong \M(N)^H\times \M(H).
\label{eq:sd}
\end{equation}

\subsection{Bogomolov multiplier} \label{sec:bog}

The following properties of $B_0(G):=\ker [H^2(G,\Q/\Z)\to
\prod_AH^2(A,\Q/\Z)]$ are taken from \cite{Bo87}, \cite{BMP}.

\begin{enumerate}

\item $B_0(G)=\ker [H^2(G,\Q/\Z)\to \prod_BH^2(B,\Q/\Z)]$, where the product is taken over
all bicyclic subgroups $B=\Z_m\times\Z_n$ of $G$  \cite{Bo87},
\cite[Cor.~2.3]{BMP}.

\item For an abelian group $A$ denote by $A_p$ its $p$-primary component.
We have $$B_0(G)=\bigoplus_pB_{0,p}(G),$$ where $B_{0,p}(G):=B_0(G)\cap \M (G)_p$.
For any Sylow $p$-subgroup $S$ of $G$ we have $B_{0,p}(G)\subseteq B_0(S)$. In particular,
if all Sylow subgroups of $G$ are abelian, $B_0(G)=0$ \cite{Bo87}, \cite[Lemma~2.6]{BMP}.

\item If $G$ is an extension of a cyclic group by an abelian group, then $B_0(G)=0$ \cite[Lemma 4.9]{Bo87}.

\item For $\gamma\in\M (G)$ consider the corresponding
central extension:
$$
1\to \Q/\Z \stackrel{i}{\to} \widetilde G_{\gamma} \to G\to 1,
$$
and denote
$$K_{\gamma }:=\{h\in \Q/\Z \,\vert\, i(h)\in \bigcap_{\chi\in\Hom (\widetilde G_{\gamma},\Q/\Z)}\ker (\chi )\}.$$
Then $\gamma$ does not belong to $B_0(G)$ if and only if some nonzero
element of $K_{\gamma}$ can be represented as a commutator of a pair
of elements of $\widetilde G_{\gamma}$ \cite[Cor.~2.4]{BMP}.

\item If $0\ne\gamma \in \M(G)$, we say that $G$ is $\gamma$-minimal if the restriction of $\gamma$ to all proper
subgroups $H\subset G$ is zero. A $\gamma$-minimal group must be a
$p$-group. We say that a $\gamma$-minimal nonabelian $p$-group  $G$
is a $\gamma$-minimal factor if for any quotient map $\rho\colon
G\to G/H$ there is no $\gamma'\in B_0(G/H)$ such that $\gamma
=\rho^*(\gamma')$ and $\gamma'$ is $G/H$-minimal. A $\gamma$-minimal
factor $G$ must be a metabelian group (i.e., $[[G,G],[G,G]]=0$) with
central series of length at most $p$, and the order of $\gamma$
in $\M(G)$ equals $p$ \cite[Theorem~4.6]{Bo87}. Moreover, if $G$ is
a $\gamma$-minimal $p$-group which is a central extension of
$G^{\text{\rm{ab}}}:=G/[G,G]$ and $G^{\text{\rm{ab}}}=(\Z_p)^n$,
then $n=2m$ and $n\ge 4$ \cite[Lemma~5.4]{Bo87}.

\end{enumerate}

\subsection{Finite simple groups} \label{sec:fin}

We need the following facts concerning finite simple groups (see,
e.g., \cite{GLS}) believing that the classification of finite simple
groups is complete.

\begin{enumerate}

\item {\it Classification.} Any finite simple group $L$ is either a group of Lie type, or an alternating group, or one of
26 sporadic groups.

\item {\it Schur multipliers.} As $L$ is perfect, it has a unique
covering group $\widetilde L$, and $L\cong \widetilde L/\M(L)$.
The Schur multipliers $\M(L)$ of all finite simple groups $L$ are
given in \cite[6.1]{GLS}.

\item {\it Automorphisms.} The group of outer automorphisms
$\Out (L):=\Aut (L)/L$ is solvable. It is abelian provided $L$ is an
alternating or a sporadic group. For groups of Lie type defined over
a finite field $F=\F_q$ the structure of $\Out (L)$ can be described
as follows. If $L$ comes from a simple algebraic group $\overline L$
defined over $\overline F$, we denote by $\overline T$ a maximal
torus in $\overline L$.
Every automorphism of $L$ is a product $idfg$ where $i$ is an inner
automorphism (identified with an element of $L$), $d$ is a diagonal
automorphism (induced by conjugation by an element $h$ of the
normalizer $N_{\overline T}(L)$, see \cite[2.5.1(b)]{GLS}), $f$ is a
field automorphism (arising from an automorphism of the field
$\overline F$), and $g$ is a graph automorphism (induced by an
automorphism of the Dynkin diagram corresponding to $L$); see
\cite[Ch.~2, \S~7]{GL} or \cite[2.5]{GLS} for more details.

The group $\Out (L)$ is a split extension of $\Outdiag (L):=\Inndiag
(L)/L$ by the group $\Phi\Gamma$, where $\Inndiag(L)$ is the group
of inner-diagonal automorphisms of $L$ (generated by all $i$'s and
$d$'s as above), $\Phi$ is the group of field automorphisms and
$\Gamma$ is the group of graph automorphisms of $L$. The group
$\cO=\Outdiag (L)$ is isomorphic to the centre of $\Lt$ by the
isomorphism preserving the action of $\Aut (L)$ and is nontrivial
only in the following cases (where $(m,n)$ stands for the greatest
common divisor of $m$ and $n$):

$L$ is of type $A_n(q)$; $\cO = \Z_{(n+1,q)}$;

$L$ is of type ${}^2A_n(q)$; $\cO = \Z_{(n+1,q-1)}$;

$L$ is of type $B_n(q)$, $C_n(q)$, or ${}^2D_{2n}(q)$; $\cO = \Z_{(2,q-1)}$;

$L$ is of type $D_{2n}(q)$; $\cO = \Z_{(2,q-1)}\times \Z_{(2,q-1)}$;

$L$ is of type ${}^2D_{2n+1}(q)$; $\cO = \Z_{(4,q-1)}$;

$L$ is of type ${}^2E_6(q)$; $\cO = \Z_{(3,q-1)}$;

$L$ is of type $E_7(q)$; $\cO = \Z_{(2,q-1)}$.

If $L$ is of type ${}^d\Sigma (q)$ for some root system $\Sigma$
($d=1,2,3)$, the group $\Phi$ is isomorphic to $\Aut (\F _{q^d})$.
If $d=1$, then $\Gamma$ is isomorphic to the group of symmetries of
the Dynkin diagram of $\Sigma$ and $\Phi\Gamma =\Phi \times \Gamma$
provided $\Sigma$ is simply-laced; otherwise, $\Gamma =1$ except if
$\Sigma =B_2, F_4$, or $G_2$ and $q$ is a power of 2, 2, or 3,
respectively, in which cases $\Phi\Gamma$ is cyclic and $[\Phi\Gamma
:\Phi]=2$. If $d\ne1$, then $\Gamma =1$.

The action of $\Phi\Gamma$ on $\cO$ is described as follows. If
$L\not\cong D_{2n}(q)$, then $\Phi$ acts on the cyclic group $\cO$
as $\Aut (\F _{q^d})$ does on the multiplicative subgroup of $\F
_{q^d}$ of the same order as $\cO$; if $L\cong D_{2n}(q)$, then
$\Phi$ centralizes $\cO$. If $L\cong A_n(q)$, $D_{2n+1}(q)$, or
$E_6(q)$, then $\Gamma =\Z_2$ acts on $\cO$ by inversion; if $L\cong
D_{2n}(q)$ and $q$ is odd, then $\Gamma$, which is isomorphic to the
symmetric group $S_3$ (for $m=2$) or to $\Z_2$ (for $m>2$) acts
faithfully on $\cO=\Z_2\times \Z_2$.

\item {\it Decorations.} It is often useful to consider groups close to finite simple groups, namely,
{\it quasisimple} and {\it almost simple} groups, as in the
statements of Theorems \ref{quasi} and \ref{almost} above. As an
example, if the simple group under consideration is $L=PSL(2,q)$,
the group $SL(2,q)$ is quasisimple and the group $PGL(2,q)$ is
almost simple. More generally, one can consider {\it semisimple}
groups (central products of quasisimple groups) and {\it nearly
simple} groups $G$, i.e., such that the generalized Fitting subgroup
$F^*(G)$ is quasisimple. $F^*(G)$ is defined as the product
$E(G)F(G)$ where $E(G)$ is the layer of $G$ (the maximal semisimple
normal subgroup of $G$) and $F(G)$ is the Fitting subgroup of $G$
(or the nilpotent radical, i.e., the maximal nilpotent normal
subgroup of $G$). The general linear group $GL(n,q)$ is an example
of a nearly simple group.

\end{enumerate}

\section{Proofs} \label{proofs}

\begin{proof}[of Theorem \ref{quasi}]
As $G$ is perfect, there exists a unique universal central covering
$\Gt$ of $G$ whose centre $Z(\Gt )$ is isomorphic to $\M (G)$ and
any other perfect central extension of $G$ is a quotient of $\Gt$.
So we can argue exactly as in \cite[Remark after Lemma 5.7]{Bo87}
and \cite{BMP}. Namely, $B_0(G)$ coincides with the collection of
classes whose restriction to any bicyclic subgroup of $G$ is zero,
see \ref{sec:bog}(1). Therefore, to establish the assertion of the
theorem, it is enough to prove that any $z\in Z(\Gt )$ can be
represented as a commutator $z=[a,b]$ of some $a,b\in \Gt$.
Moreover, it is enough to prove that such a representation exists
for all elements $z$ of prime power order, see \ref{sec:bog}(2).

It remains to apply the results of Blau \cite{Bl} who classified all
elements $z$ having a fixed point in the natural action on the set
of conjugacy classes of $\Gt$ (such elements evidently admit a
needed representation as a commutator):

\begin{thm} {\text{\rm(\cite[Theorem 1]{Bl})}} \label{th:Blau}
Assume that $G$ is a quasisimple group and let $z\in Z(G)$. Then one
of the following holds:

(i) $\order(z)=6$ and $G/Z(G)\cong A_6, A_7, Fi_{22}, PSU(6,2^2),$ or ${}^2E_6(2^2)$;

(ii) $\order(z)=6$ or $12$ and $G/Z(G)\cong PSL(3,4), PSU(4,3^2)$ or $M_{22}$;

(iii) $\order(z)=2$ or $4$, $G/Z(G)\cong PSL(3,4)$, and $Z(G)$ is noncyclic;

(iv) there exists a conjugacy class $C$ of $G$ such that $Cz=C$.
\end{thm}

This theorem implies that the only possibility for an element of $G$
of prime power order to act on the set of conjugacy classes without
fixed points is case (iii) where $\Gt /Z(\Gt )\cong PSL(3,4)$ and
$z$ is an element of order 2 or 4. So the classes $\gamma\in
H^2(G,\Q /\Z)$ corresponding to such $z$'s are the only candidates
for nonzero elements of $B_0(G)$.

A more detailed analysis of the case $PSL(3,4)$, where $Z:=Z(\Gt
)\cong \Z_3\times\Z_4\times\Z_4$, is sketched in
\cite[Remark (2) after Theorem 1]{Bl}. The result (rechecked by
MAGMA computations) looks as follows: all elements of $Z$ of orders 2 and 3
fix some conjugacy class of $\Gt$, all elements of orders 6 and 12
act without fixed points, and of the twelve elements of order 4
exactly six fix a conjugacy class of $\Gt$.

First note that this description implies $B_0(PSL(3,4))=0$. Indeed,
the criterion given in \ref{sec:bog}(4) can be rephrased for a
quasisimple group $G$ as follows: $B_0(G)=0$ if and only if $Z(\Gt
)$ has a system of generators each of those can be represented as a
commutator of a pair of elements of $\Gt$. It remains to notice that
if $G=PSL(3,4)$, each 5-tuple of elements of order 4 in $Z$
generates $\Z_4\times \Z_4$ because the subgroup of the shape
$\Z_4\times \Z_2$ contains only 4 elements of order 4 (I am indebted
to O.~Gabber for this argument). (Another way to prove that
$B_0(PSL(3,4))=0$ was demonstrated in \cite{BMP} where Lemma 5.3
establishes a stronger result: vanishing of $B_0(S)$, where $S$ is a
2-Sylow subgroup of $PSL(3,4)$.)

The only cases where the condition of the above mentioned criterion
breaks down are those where $G$ is a 12- or 4-cover of $PSL(3,4)$.
Indeed, in these cases we have $G=\Gt /Z$ where $Z$ is generated
either by an element of order 4 (which may be not representable
as a commutator) or of order 12 (which cannot be representable
as a commutator). Thus in these cases we have $B_0(G)\neq 0$.
More precisely, in both cases we have $B_0(G)=\Z_2$ because the subgroup
generated by commutators is of index 2 in $G$ (I thank O.~Gabber for
this remark).
\end{proof}

\begin{rem} \label{Th}
It is interesting to compare \cite[Lemma 3.1]{BMP} with a theorem
from the PhD thesis of Robert Thompson \cite[Theorem 1]{Th}.
\end{rem}

\begin{proof}[of Theorem \ref{almost}]
Let $L\subseteq G\subseteq \Aut (L)$ where $L$ is a simple group.
Clearly, it is enough to prove the theorem for $G=\Aut (L)$. The
group $\Out (L)=\Aut(L)/L$ of outer automorphisms of $L$ acts on $\M
(L)$, and since $L$ is perfect, we have an isomorphism
\begin{equation}
\M (G)\cong \M (L)^{\Out (L)} \times \M (\Out (L))\label{eq:semi}
\end{equation}
(see \eqref{eq:sd}).

\begin{lemm} \label{out}
$B_0(\Out (L))=0$.
\end{lemm}

\begin{proof}[of Lemma \ref{out}]
We maintain the notation of Section \ref{sec:fin}. If $\Out (L)$ is
abelian, the statement is obvious. This includes the cases where $L$
is an alternating or a sporadic group. So we may assume $L$ is of
Lie type. If $\cO=1$, i.e., $L$ is of type $E_8$, $F_4$, or $G_2$,
the result follows immediately. If the group $\Phi\Gamma$ is cyclic,
the result follows because $\cO$ is abelian (see \ref{sec:fin}(3)).
This is the case for all groups having no graph automorphisms, in
particular, for all groups of type $B_n$ or $C_n$ ($n\ge 3)$, $E_7$,
and for all twisted forms. For the groups of type $B_2$, the group
$\Phi\Gamma$ is always cyclic. It remains to consider the cases
$A_n$, $D_n$, and $E_6$. In the case $L=E_6$ all Sylow $p$-subgroups
of $\Out (L)$ are abelian, and the result holds. Let $L=D_{2m}(q)$.
If $q$ is even, we have $\cO =1$, $\Gamma=\Z_2$ (if $m>2$) or $S_3$
(if $m=2$); in both cases the Sylow $p$-subgroups of $\Out (L)$ are
abelian, and we are done. If $q$ is odd, we have $\cO =\Z_2\times
\Z_2$, and $\Phi$ centralizes $\cO$ (see Section \ref{sec:fin}), so
every Sylow $p$-subgroup of $\Out (L)$ can be represented as an
extension of a cyclic group by an abelian group, and we conclude as
above. Finally, let $L$ be of type $A_n(q)$ or $D_{2m+1}(q)$. Then
we have $\cO=\Z_h$, $h=(n+1,q-1)$ or $h=(4,q-1)$, respectively,
$\Gamma=\Z_2$, $\Phi = \Aut (\F _q)$. The action of both $\Gamma$
and $\Phi$ on $\cO$ may be nontrivial: $\Gamma$ acts by inversion,
$\Phi$ acts on $\cO$ as $\Aut (\F _q)$ does on the multiplicative
subgroup of $\F_q$ of the same order as $\cO$. Hence we can
represent the metabelian group $\Out (L)$ in the form
\begin{equation}
1\to V\to \Out (L)\to A\to 1,
\label{eq:meta}
\end{equation}
where $V$, the derived subgroup of $\Out (L)$, is isomorphic to a
cyclic subgroup $\Z_c$ of $\cO$, and the abelian quotient $A$ is of
the form $\Z_a\times \Z_b \times \Z_2$ for some integers $a,b,c$.
Since it is enough to establish the result for a Sylow 2-subgroup,
we may assume that $a$, $b$ and $c$ are powers of $2$. Then the
statement of the lemma follows from the properties of
$\gamma$-minimal elements described in Section \ref{sec:bog}.
Indeed, if $\gamma$ is a nonzero element of $B_0(G)$ and $G$ is
$\gamma$-minimal, then $G$ is metabelian, both $V$ and $A$ are of
exponent $p$, and in any representation of $G$ in the form
(\ref{eq:meta}) the group $A$ must have an even number $s=2t$ of direct
summands $\Z_p$ with $t\ge 2$. However, if $G$ is a Sylow
2-subgroup of $\Out (L)$, this is impossible because $A$ contains
only three direct summands. Thus $B_0(Syl_2(\Out (L)))=0$, and so
$B_0(\Out (L))=0$. The lemma is proved.
\end{proof}

We can now finish the proof of the theorem. Let $\gamma$ be a
nonzero element of $B_0(G)$. Using the isomorphism (\ref{eq:semi}),
we can represent $\gamma$ as a pair $(\gamma_1, \gamma_2)$ where
$\gamma_1\in \M(L)$, $\gamma_2\in \M(\Out (L))$. Restricting to the
bicyclic subgroups of $G$, we see that $\gamma_1\in B_0(L)$,
$\gamma_2\in B_0(\Out (L))$, and the result follows from Theorem
\ref{quasi} and Lemma \ref{out}.
\end{proof}


\noindent {\bf Acknowledgements}. The author's research was
supported in part by the Minerva Foundation through the Emmy Noether
Research Institute of Mathematics, the Israel Academy of Sciences
grant 1178/06, and a grant from the Ministry of Science, Culture
and Sport (Israel) and the Russian Foundation for Basic Research
(the Russian Federation). This paper was mainly written during the visit to
the MPIM (Bonn) in August--September 2007 and completed during the visit
to ENS (Paris) in April--May 2008. The support of these institutions is highly
appreciated.

My special thanks are due to O.~Gabber who noticed a gap in an earlier version
of the paper and helped to fill it. I am also grateful to M.~Conder, D.~Holt
and A.~Hulpke for providing representations of $PSL(3,4)$ necessary for
MAGMA computations, and to N. A. Vavilov for useful correspondence.


\begin{thebibliography}{Bo07}

\bibitem[Bl]{Bl}
{\sc H. I. Blau}, {\it A fixed-point theorem for central elements in quasisimple groups},
Proc. Amer. Math. Soc. {\bf 122} (1994) 79--84.

\bibitem[Bo87]{Bo87}
{\sc F. A. Bogomolov}, {\it The Brauer group of quotient spaces by linear group actions},
Izv. Akad. Nauk. SSSR Ser. Mat. {\bf 51} (1987) 485--516; English transl. in
Math. USSR Izv. {\bf 30} (1988) 455--485.


\bibitem[Bo92]{Bo92}
{\sc F. A. Bogomolov}, {\it Stable cohomology of groups and algebraic varieties},
Mat. Sb. {\bf 183} (1992) 1--28; English transl. in
Sb. Math. {\bf 76} (1993) 1--21.


\bibitem[Bo07]{Bo07}
{\sc F. Bogomolov}, {\it Stable cohomology of finite and profinite groups},
``Algebraic Groups'' (Y.~Tschinkel, ed.), Universit\"atsverlag G\"ottingen, 2007, 19--49.


\bibitem[BMP]{BMP}
{\sc F. Bogomolov, J. Maciel, T. Petrov},
{\it Unramified Brauer groups of finite simple groups of Lie type $A_{\ell}$},
Amer. J. Math. {\bf 126} (2004) 935--949.


\bibitem[CTS]{CTS} {\sc J.-L. Colliot-Th\'el\`ene, J.-J. Sansuc},
{\it The rationality problem for fields of
invariants under linear algebraic groups (with special regards to the Brauer group)},
Proc. Intern. Colloquium on Algebraic Groups and Homogeneous Spaces (Mumbai
2004)
(V.~Mehta, ed.), TIFR Mumbai, Narosa Publishing House, 2007, 113--186.

\bibitem[GS]{GS}
{\sc P. Gille, T. Szamuely}, {\it Central Simple Algebras and Galois
Cogomology}, Cambridge Univ. Press, Cambridge, 2006.


\bibitem[GL]{GL}
{\sc D. Gorenstein, R. Lyons}, {\it The Local Structure of Finite Groups
of Characteristic $2$ Type}, Mem. Amer. Math. Soc., vol.~276,
Providence, RI, 1983.


\bibitem[GLS]{GLS}
{\sc D. Gorenstein, R. Lyons, R. Solomon}, {\it The Classification
of the Finite Simple Groups}, Number 3, Math. Surveys and
Monographs, vol.~40, no.~3, Amer. Math. Soc., Providence, RI,
1998.

\bibitem[Ka]{Ka}
{\sc G. Karpilovsky},
{\it Group Representations}, vol.~2,
North-Holland Math. Studies {\bf 177}, North-Holland, Amsterdam et al., 1993.

\bibitem[Sa]{Sa}
{\sc D. J. Saltman},
{\it Noether's problem over an algebraically closed field},
Invent. Math. {\bf 77} (1984) 71--84.

\bibitem[Sh]{Sh}
{\sc I. R. Shafarevich}, {\it The L\"uroth problem},
Trudy Mat. Inst. Steklov {\bf 183} (1990) 199--204; English transl.
Proc. Steklov Inst. Math. {\bf 183} (1991) 241--246.

\bibitem[Th]{Th}
{\sc R. C. Thompson},
{\it Commutators in the special and general linear groups},
Trans. Amer. Math. Soc. {\bf 101} (1961) 16--33.



\end{thebibliography}
\end{document}